\documentclass{article}
\usepackage{amsmath}[1996/11/01]
\usepackage{amssymb,amsthm,amsxtra}

\def\[#1\]{\begin{equation}#1\end{equation}}
\makeatletter
\def\beq{%
   \relax\ifmmode
      \@badmath
   \else
      \ifvmode
         \nointerlineskip
         \makebox[.6\linewidth]%
      \fi
      $$
   \fi
}
\def\eeq{%
   \relax\ifmmode
      \ifinner
         \@badmath
      \else
         $$
      \fi
   \else
      \@badmath
   \fi
   \ignorespaces
}

\def\enddisplaymath{\eeq\global\@ignoretrue}
\makeatother

\newtheorem{thm}{Theorem}
\newtheorem*{thmi}{Theorem}
\newtheorem{cor}[thm]{Corollary}
\newtheorem{lem}[thm]{Lemma}

\theoremstyle{remark}
\newtheorem*{rem}{Remark}
\newtheorem{rems}{Remark}[thm]

\theoremstyle{definition}

\numberwithin{equation}{section}
\numberwithin{thm}{section}

\renewcommand{\Re}{\operatorname{Re}}

\newcommand{\Z}{\mathbb Z}
\newcommand{\Q}{\mathbb Q}

\begin{document}

\title{{\bf New asymptotic bounds for self-dual codes and lattices}}
\author{Eric M. Rains\\AT\&T Research\\rains@research.att.com}
\date{April 12, 2001}
\maketitle

\begin{abstract}
We give an independent proof of the Krasikov-Litsyn bound $d/n\lesssim
(1-5^{-1/4})/2$ on doubly-even self-dual binary codes.  The technique used
(a refinement of the Mallows-Odlyzko-Sloane approach) extends easily to
other families of self-dual codes, modular lattices, and quantum codes; in
particular, we show that the Krasikov-Litsyn bound applies to singly-even
binary codes, and obtain an analogous bound for unimodular lattices.  We
also show that in each case, our bound differs from the true optimum by an
amount growing faster than $O(\sqrt{n})$.

\end{abstract}

\section{Introduction}

In \cite{MallowsCL/SloaneNJA:1973},
\cite{MallowsCL/OdlyzkoAM/SloaneNJA:1975}, Mallows, Odlyzko, and Sloane
proved the following result:

\begin{thmi} Let $C$ be a doubly-even binary self-dual code of length $n$
and minimum distance $d$.  Then $d\le 4[n/24]+4$.  Morever, for any
constant $b$, one has $d\le n/6-b$ for sufficiently large $n$.
\end{thmi}
\noindent as well as analogous results for ternary codes and even unimodular
lattices.  The first claim has since been extended to singly-even binary
self-dual codes \cite{me:G}, and analogous results have been obtained for
even strongly modular lattices \cite{QuebbemannH-G:1995},
\cite{QuebbemannH-G:1997} and odd strongly modular lattices \cite{me:M},
including the odd unimodular case.

Regarding the asymptotic claim, essentially the only improvement is the
bound of Krasikov and Litsyn \cite{KrasikovI/LitsynS:2000}:

\begin{thmi} Let $C_i$ be a family of doubly-even binary self-dual codes
of length tending to infinity.  Then
\[
\limsup_{i\to\infty} \frac{d(C_i)}{n(C_i)} \le \frac{1-5^{-1/4}}{2}.
\]
\end{thmi}
\noindent However, it is unclear to what extent their argument extends to
the other cases of interest (especially the lattice cases).

In the present paper, we give a new technique for deriving bounds on
self-dual codes and modular lattices.  In the doubly-even binary case, our
bound is precisely the Krasikov-Litsyn bound; the difference is that our
technique easily generalizes.  The basic idea (following
\cite{MallowsCL/OdlyzkoAM/SloaneNJA:1975}) is to use invariant theory to
construct linear relations that must be satisfied by the weight enumerator
of a self-dual code.  In the doubly-even binary case, the simplest such
relation gives $d\le 4[n/24]+4$; to obtain their asymptotic improvement,
Mallows, Odlyzko, and Sloane also take into account the second-simplest
relation.  Somewhat surprisingly, our bounds also only use these two
relations; equivalently, we only consider the coefficients of the weight
enumerator up to weight $4[n/24]+8$.

To compute the coefficients of these relations, we use the
B\"urmann-Lagrange theorem to express them as coefficients of certain power
series, which we asymptotically analyze via Cauchy's integral and the
saddle-point method.  Under suitable conditions on the power series, we can
then show that their coefficients in certain ranges are asymptotically
uniformly positive, and thus give a contradiction unless the minimum
distance bound holds.  This necessitates a certain amount of analysis,
which we deal with in Section \ref{sec:lemmas}.  Then, in Section
\ref{sec:main}, we give our main theorem, Theorem \ref{thm:main1}.  This
is stated in some generality (regarding the minimum valuations of certain
families of power series with nonnegative coefficients), so as to include
most of our applications as special cases.  We state these special cases in
Section \ref{sec:apps1}; see Theorem \ref{thm:codes2} for (most) codes,
Theorem \ref{thm:latts2} for even modular lattices, and Theorems
\ref{thm:Z4HL} and \ref{thm:Z4E2} for codes over $\Z_4$.  The remaining
applications not directly dealt with by Theorem \ref{thm:main1} are
considered in Section \ref{sec:apps2}; see Theorem
\ref{thm:codes1} for (singly-even) self-dual binary codes,
Theorem \ref{thm:latts1} for odd modular lattices, and Theorems
\ref{thm:quantq} and \ref{thm:quant2} for quantum codes.

As we remarked, our improvements on the main Mallows-Odlyzko-Sloane bounds
are obtained by considering the first two relations coming from invariant
theory, rather than just the first.  In Section \ref{sec:hermite}, we
consider the possibility of improving the bounds by using the first $k$
relations.  In fact, it turns out that, despite the significant improvement
between $k=1$ and $k=2$, increasing $k$ beyond $2$ does not give a better
bound on $\limsup d/n$.  We do, however, obtain a slight lower-order
improvement (Theorem \ref{thm:hermite}); increasing $k$ gives an
$O(n^{-1/2})$ improvement on the bound on $d/n$.  In particular, each of
our bounds differs from the true optimal minimum distance (norm) by an
amount growing faster than $O(\sqrt{n})$.

{\bf Acknowledgements}.  The author would like to thank H. Landau,
A. M. Odlyzko, and N. J. A. Sloane for helpful discussions regarding
Section \ref{sec:lemmas}, especially Lemma \ref{lem:mainlem}, as well
as I. Duursma for pointing out that Krasikov and Litsyn had improved
their earlier bound to the one stated above.

\section{Lemmas}\label{sec:lemmas}

We will use the notation $[t^j] f(t)$ to refer to the coefficient of
$t^j$ in the (formal) Laurent series $f(t)$.  We also use the notation
\[
f\succeq g
\]
to say that $[t^j] (f(t)-g(t))\ge 0$ for all $j$.

\begin{lem}
Let $f(t)$ be a Laurent series convergent on an annulus $0\le
r_1<|t|<r_2\le \infty$.  If $f\succeq 0$, then we have the bound
\[
|f(t)|\le f(|t|)
\]
valid on the annulus.  If $[t^j] f(t)$ and $[t^{j+1}] f(t)$ are both
nonzero for some $j$, then equality can hold only when $t=|t|$.
\end{lem}

\begin{proof}
We have:
\[
|f(t)| = |\sum_j f_j t^j|\le \sum_j f_j |t|^j = f(|t|),
\]
with equality only when there exists $\alpha$ such that $t^j=\alpha |t|^j$
for all $j$ with $f_j$ nonzero.  Dividing two consecutive such equations,
we obtain $t=|t|$ as required.
\end{proof}

The Hadamard three-circles theorem then immediately implies that
$\log(f(e^s))$ is strictly convex, for $r_1<e^s<r_2$.  In fact, we have the
slightly stronger statement:

\begin{lem}\label{lem:logconvex}
let $f(t)$ be a Laurent series convergent on an annulus $0\le
r_1<|t|<r_2\le \infty$, and not proportional to $z^n$ for any $n$.
If $f\succeq 0$, then for $r_1<e^s<r_2$,
\[
\frac{d^2}{ds^2} \log(f(e^s)) > 0.
\]
\end{lem}

\begin{proof}
Setting $r=e^s$, we have:
\[
\frac{d^2}{ds^2} \log(f(e^s))
=
\left(r \frac{d}{dr}\right)^2 \log(f(r))
=
\frac{\left(r\frac{d}{dr}\right)^2 f(r)}{f(r)}
-
\left(\frac{r\frac{d}{dr} f(r)}{f(r)}\right).
\]
Now, $f(r)$ is positive on $r_1<r<r_2$, so we may freely multiply by
$f(r)^2$; we thus need to show that
\[
f(r) \left(r\frac{d}{dr}\right)^2 f(r)
-
\left(r\frac{d}{dr} f(r)\right)^2
>
0.
\]
Now, the left-hand-side has a Laurent series convergent in the original
annulus, namely
\begin{align}
(\sum_j f_j t^j)(\sum_k k^2 f_k t^k)
-
(\sum_j j f_j t^j)(\sum_k k f_k t^k)
&=
\sum_j \sum_k k^2 f_j f_k t^{j+k}
-
\sum_j \sum_k jk f_j f_k t^{j+k}\\
&=
\sum_j \sum_k \frac{(j-k)^2}{2} f_j f_k t^{j+k}.
\end{align}
Since this $\succeq 0$, and has at least one nonzero coefficient, the
desired inequality follows.
\end{proof}


\begin{lem}\label{lem:mainlem}
Let $F,G$ be real power series both convergent in the circle
$|t|<r_0$.  Suppose furthermore that $G\succeq 0$, with both $G(0)$ and
$G'(0)$ positive.  Then for any compact subset $I\subset [0,r_0)$ on which $F$
is positive, and for all sufficiently large $n$ (depending on $I$):
\[
[t^m] F(t) G(t)^n > 0
\]
whenever $m/n\in S(I)$, with $S(t):=t G'(t)/G(t)$.
\end{lem}

\begin{proof}
It will suffice to show that
\[
[t^{S(r)n}] F(t) G(t)^n = (F(r)+o(1)) [t^{S(r)n}] G(t)^n,
\]
with error uniform on any interval $r\in [0,a]$ with $a<r_0$, since then on
$I$, $F(r)$ is bounded away from 0, while the error converges uniformly to
0.  We split into two cases: $0\le r\le n^{-2/3}$ and $n^{-2/3}\le r\le a$.

In the first region, we claim that for all sufficiently large $n$, and for
$0\le r\le n^{-2/3}$,
\[
[t^{S(r)n}] F(t) G(t)^n
=
(F(0)+O(n^{-1/3})) G(0)^n \frac{(n S'(0))^{S(r)n}}{(S(r)n)!}
\]
with error uniform in $r$.  By Cauchy's residue theorem,
\[
G(0)^{-n} [t^{S(r)n}] F(t) G(t)^n
=
\frac{1}{2\pi i}
\int_{|t|= S(r)/S'(0)} t^{-S(r)n} F(t) \left(\frac{G(t)}{G(0)}\right)^n.
\]
Now, $|t|=O(n^{-2/3})$, so we have the uniform estimates
\begin{align}
\log G(t)-\log G(0) &= S'(0) t + O(n^{-4/3}),\\
F(t) &= F(0) + O(n^{-2/3}),
\end{align}
and thus
\[
F(t) \left(\frac{G(t)}{G(0)}\right)^n = e^{S'(0) n t} (F(0)+O(n^{-1/3})).
\]

Now,
\[
\frac{1}{2\pi i} \int_{|t|=S(r)/S'(0)} t^{-S(r)n} F(0) e^{S'(0) n t} dt/t
=
F(0) \frac{(n S'(0))^{S(r)n}}{(S(r)n)!},
\]
so it remains to show that
\[
(S(r) n)! (S'(0) n)^{-S(r)n}
\frac{1}{2\pi i} \int_{|t| = S(r)/S'(0)} |\exp(S'(0) n t)| \frac{dt}{t}
\]
is bounded.  But, setting $m=S(r)n$ and rescaling $t$, this is
\[
\frac{m!}{m^m}
\frac{1}{2\pi i} \int_{|t|=1} e^{m\Re(t)} \frac{dt}{t}
=
\frac{m! I_0(m)}{m^m} = 1+O(1/m),
\]
by the known asymptotics of Bessel functions.

\bigskip

We now consider the case $n^{-2/3}\le r\le a$.  Here, we claim
\[
[t^{S(r)n}] F(t) G(t)^n
=
(2\pi n r S'(r))^{-1/2} r^{-S(r)n} G(r)^n
(F(r)+O(n^{-1/60})),
\]
again with uniform error.  Again, Cauchy's integral gives
\[
r^{S(r)n} G(r)^{-n} [t^{S(r)n}] F(t) G(t)^n
=
\frac{1}{2\pi}
\int_{-\pi}^\pi
e^{-i S(r) n \theta}
F(r e^{i\theta})
\left(\frac{G(r e^{i\theta})}{G(r)}\right)^n
d\theta.
\]

Since 
\[
\frac{\log |G(re^{i\theta})|-\log G(r)}{r\theta^2}
\]
is continuous and negative for $r\in [0,b]$, $\theta\in [-\pi,\pi]$,
there exists a positive constant $C$ such that
\[
\log |G(re^{i\theta})|-\log G(r) \le -C r\theta^2
\]
in that region.

In particular, when $|\theta|\ge n^{-9/20} r^{-1/2}$, we find
$n r \theta^2 \ge n^{1/10}$, and thus
\[
\left|\frac{G(re^{i\theta})}{G(r)}\right|^n\le e^{-C n^{1/10}}.
\]
Since $|F(re^{i\theta})|$ is bounded, we find that the contribution of
this region to the integral is negligible.

Now, consider the region $|\theta|\le n^{-9/20} r^{-1/2}$.
Since
\[
\frac{\log G(re^{i\theta})-\log G(r)-iS(r)\theta+r S'(r)\frac{\theta^2}{2}}
{r\theta^3}
\]
is continuous for $r\in [0,b]$, $\theta\in [-\pi,\pi]$, we find the
uniform estimate
\[
\log G(re^{i\theta})-\log G(r)-iS(r)\theta+\frac{r S'(r)\theta^2}{2}
=
O(r\theta^3).
\]
or upon exponentiation,
\[
e^{-iS(r)n\theta} G(r)^{-n} G(re^{i\theta})^n
=
e^{-nrS'(r)\theta^2/2}
(1+O(n r\theta^3))
=
e^{-nrS'(r)\theta^2/2} (1+O(n^{-1/60})),
\]
Similarly,
\[
F(r e^{i\theta}) = F(r) + O(n^{-7/60}),
\]
and thus
\[
F(r e^{i\theta}) e^{-iS(r)n\theta} G(r)^{-n} G(re^{i\theta})^n
=
e^{-n rS'(r)\theta^2/2} (F(r)+O(n^{-1/60})).
\]
Since
\[
\frac{1}{2\pi} \int_{|\theta|\le n^{-9/20} r^{-1/2}}
e^{-n r S'(r)\theta^2/2} d\theta
=
(2\pi n r S'(r))^{-1/2}
+
O(e^{-n^{1/10} S'(r)/2}),
\]
the claim follows.
\end{proof}

\begin{rem}
Note that we only used the fact $G\succeq 0$ through the conclusions of the
previous two lemmas.
\end{rem}

Away from $0$, we can give much stronger estimates (which will be used in
the final section):

\begin{lem}\label{lem:asymptseries}
Fix radii $0<r_1<r_2<\infty$, and let $F(t)$ and $G(t)$ be Laurent series
convergent on a neighborhood of the annulus $r_1\le |t|\le r_2$; suppose
further that $G(t)\succeq 0$ and has two consecutive nonzero coefficients.
Let $S(r):=r G'(r)/G(r)$, and for each $r_1\le r\le r_2$, define a power
series $\gamma(r,x)$ by $\gamma(r,0)=0$, $\gamma_x(r,0)>0$ and
\[
\log G(r e^{i \gamma(r,x)}) = \log G(r) + i S(r) \gamma(r,x)
-
r S'(r) \frac{x^2}{2}.
\]
Then for all integers $k>0$ and for $r_1\le r\le r_2$, we have the
asymptotic estimate
\[
\sqrt{2\pi r S'(r) n} r^{S(r) n} G(r)^{-n} [t^{S(r) n}] F(t) G(t)^n
=
\sum_{0\le j\le k-1}
\frac{(r S'(r) n)^{-j}}{2^j j!}
\left(\frac{d^{2j}}{dx^{2j}} \gamma_x(r,x) F(r e^{i\gamma(r,x)})\right)_{x=0}
+
O(n^{-k}),
\]
with error uniform in $r$.
\end{lem}

\begin{proof}
Note that
\[
\log G(r e^{i y}) = \log G(r) + i S(r) y - r S'(r) \frac{y^2}{2} + O(y^3);
\]
since $r S'(r)>0$, we conclude that $\gamma(r,x)$ converges for
$|x|\le x_0$ for some $x_0>0$ independent of $r$, and satisfies
$\gamma_x(r,0)=1$.

By Laurent's theorem,
\[
\sqrt{2\pi r S'(r) n} r^{S(r) n} G(r)^{-n} [t^{S(r) n}] F(t) G(t)^n
=
\sqrt{\frac{r S'(r) n}{2\pi}}
\int_{-\pi}^\pi
F(r e^{i\theta})
\left(\frac{e^{-i S(r)\theta} G(r e^{i\theta})}{G(r)}\right)^n
d\theta.
\]
Now, as before, we can restrict the integral to any uniform neighborhood of
0, with exponentially small error.  In particular, we may restrict
to a neighborhood $|\theta|\le \theta_0$ affording the change of variable
$\theta = \gamma(r,x)$.  The integral thus becomes
\[
\frac{r S'(r) n}{2\pi}
\int_{-x_0}^{x_0}
\gamma_x(r,x) F(r e^{i \gamma(r,x)}) e^{-n r S'(r) x^2/2} dx.
\]
Now, we have the uniform estimate
\[
\gamma_x(r,x) F(r e^{i \gamma(r,x)})
=
\sum_{0\le j\le 2k-1}
\frac{x^j}{j!} (\frac{d^j}{dx^j} \gamma_x(r,x) F(r e^{i \gamma(r,x)}))_{x=0}
+
O(x^{2k});
\]
since
\[
\sqrt{r S'(r) n}{2\pi}
\int_{-x_0}^{x_0}
x^{2k} e^{-n r S'(r) x^2/2} dx
\le
\sqrt{r S'(r) n}{2\pi}
\int_{-\infty}^\infty x^{2k} e^{-n r S'(r) x^2/2} dx
=
(r S'(r) n)^{-k} \frac{(2k)!}{2^k k!},
\]
the contribution of the error term is as required.  Once we remove
this term, the integral can be extended to $\infty$, again giving uniform
exponentially small error.  Evaluating the resulting Gaussian integral gives
the desired result.
\end{proof}

For our purposes, we will need a version of this valid in the neighborhood
of a zero of $F$.  Define polynomials $h_k(x)$ for integers $k\ge 0$ via
the generating function
\[
\sum_{k\ge 0} h_k(x) \frac{t^k}{k!} = e^{tx-t^2/2};
\]
in particular, $h_k(x)$ is a monic polynomial of degree $k$, and
\begin{align}
h_k(x) &= \sum_{0\le j\le [k/2]}
(-1)^j
\frac{k!}{2^j j!}
\frac{x^{k-2j}}{(k-2j)!}\\
&=
\sum_{0\le j\le k} \binom{k}{j} x^{k-j} h_j(0).
\end{align}
(Thus $h_k(x)$ are just rescaled Hermite polynomials.)  In terms of
$h_j(0)$, the above estimate becomes (replacing $k$ by $k/2$):
\begin{multline}
\sqrt{2\pi r S'(r) n} r^{S(r) n} G(r)^{-n} [t^{S(r) n}] F(t) G(t)^n
=\\
\sum_{0\le j\le k}
\frac{(r S'(r) n)^{-j/2} i^{-j} h_j(0)}{j!}
\left(\frac{d^j}{dx^j} \gamma_x(r,x) F(r e^{i\gamma(r,x)})\right)_{x=0}
+
O(n^{-(k+1)/2}).
\label{eq:asympt_series2}
\end{multline}

\begin{cor}\label{cor:asympt_hermite}
With hypotheses as in Lemma \ref{lem:asymptseries}, suppose $F$ has a zero
of order $k$ at the point $r_0\in [r_1,r_2]$.  Then we have the uniform
asymptotic estimate
\[
\frac{[t^{S(r) n}] F(t) G(t)^n}{[t^{S(r) n}] G(t)^n}
=
\frac{F^{(k)}(r_0)}{k!}
\frac{h_k\bigl((r-r_0)\sqrt{S'(r_0) n/r_0}\bigr)}{(S'(r_0) n/r_0)^{k/2}}
+
O(\max(|r-r_0|,n^{-1/2})^{k+1}),
\]
valid for $r_1\le r\le r_2$.
\end{cor}

\begin{proof}
Since the main term of the estimate has order $O(|r-r_0|^k)$, we can
tolerate a multiplicative error of order $1+O(n^{-1/2})$.  In particular,
we may replace $[t^{S(r) n}] G(t)^n$ by its first-order estimate.
We thus need to estimate \eqref{eq:asympt_series2}.

Now,
\begin{align}
\gamma_x(r,x) F(r e^{i\gamma(r,x)})
&=
\sum_{0\le l}
((r \frac{d}{dr})^l F)(r) \gamma_x(r,x) \gamma(r,x)^l \frac{i^l}{l!}\\
&=
\sum_{0\le l}
((r \frac{d}{dr})^l F)(r) \frac{d}{dx} \gamma(r,x)^{l+1} \frac{i^l}{(l+1)!}\\
&=
\sum_{0\le j}
\frac{x^j}{j!}
\sum_{0\le l\le j}
\frac{i^l}{(l+1)!}
((r \frac{d}{dr})^l F)(r) \left((\frac{d}{dx})^{j+1} h(r,x)^{l+1}\right)_{x=0}.
\end{align}
Thus we obtain
\[
\sum_{0\le j\le k}
\frac{(r S'(r) n)^{-j/2} h_j(0) i^{-j}}{j!}
\sum_{0\le l\le j}
\frac{i^l}{(l+1)!}
((r \frac{d}{dr})^l F)(r) \left((\frac{d}{dx})^{j+1} \gamma(r,x)^{l+1}\right)_{x=0}.
\]

Now, for $0\le l\le k$,
\[
(r \frac{d}{dr})^l F(r) = O((r-r_0)^{k-l}),
\]
and for $l\le j$, we have:
\[
\left((\frac{d}{dx})^{j+1} \gamma(r,x)^{l+1}\right)_{x=0} = O(1).
\]
In particular, the $j,l$ term gives a contribution of order
$O(|r-r_0|^{k-l} n^{-j/2})$.
For fixed $l$, the contributions get smaller as $j$ increases.
We thus find that the terms with $l<j\le k$ are dominated by
the terms with $j=l+1$, of order
\[
O(|r-r_0|^{k-l} n^{-(l+1)/2}) = O(\max(|r-r_0|,n^{-1/2})^{k+1}).
\]
It remains to consider the terms with $j=l$, that is,
\[
\sum_{0\le l\le k}
\frac{(r S'(r) n)^{-l/2} h_l(0)}{l!}
((r \frac{d}{dr})^l F)(r).
\]
If we replace $r \frac{d}{dr}$ by $r_0\frac{d}{dr}$, $F(r)$ by
$F^{(k)}(r_0) (r-r_0)^k/k!$ and $r S'(r)$ by $r_0 S'(r_0)$, the resulting
error is again $O(\max(|r-r_0|,n^{-1/2})^{k+1})$.  We thus obtain
\[
\frac{F^{(k)}(r_0)}{k!}
\sum_{0\le l\le k}
\sqrt{S'(r_0) n/r_0})^{-l}
\frac{h_l(0)}{(S'(r_0) n/r_0)^{l/2}}
\binom{k}{l} (r-r_0)^{k-l}
=
\frac{F^{(k)}(r_0)}{k!}
\frac{h_k\bigl(\sqrt{S'(r_0) n/r_0} (r-r_0)\bigr)}{(S'(r_0) n/r_0)^{k/2}}
\]
as required.
\end{proof}

We conclude with one more analytical lemma, used in Section \ref{sec:apps2}
below.

\begin{lem}\label{lem:shadow}
Let $F(t)$ and $G(t)$ be real power series convergent on the circle
$|t|<r$.  Suppose $G\succeq 0$, with $G(0)$, $G'(0)$ both positive.
Let $I$ be a compact subset of $[0,r)$ on which $F$ is positive.
Then for all sufficiently large $l$, the holomorphic functions
\[
F(a) G(a)^l + F(-a) G(-a)^l\quad\text{and}\quad
\frac{F(a) G(a)^l - F(-a) G(-a)^l}{a}
\]
are positive on $I$.
\end{lem}

\begin{proof}
We first note that if $I$ does not contain 0, then $|G(-a)/G(a)|$ is
bounded below 1 on $I$, while $|F(-a)/F(a)|$ is bounded; the result follows
immediately.  If $I$ is the single point $\{0\}$, then $F(0)>0$ by
assumption, so $F(x)>0$ for sufficiently small positive $x$; we may thus
enlarge $I$ while maintaining the hypotheses.  We may thus take $I$ of the
form $[0,\epsilon]$ with $\epsilon>0$; moreover, if the theorem is true for
$[0,\epsilon']$ with $0<\epsilon'<\epsilon$, it is true for $[0,\epsilon]$.

Since $F(0)$ and $G(0)$ are both positive, we may choose $\epsilon$
so that $F(a)$ and $G(a)$ are positive on $[-\epsilon,\epsilon]$.
In particular, the first function is thus positive on $[0,\epsilon]$ for
all $l$; it remains to consider the second function.

Now, clearly $(F(a)-F(-a))/a$ is analytic on $[-\epsilon,\epsilon]$,
while
\[
\frac{\log G(a)-\log G(-a)}{a}
\]
is real analytic and positive on $[-\epsilon,\epsilon]$.  In particular,
there exist positive constants $C_1$ and $C_2$ such that
\begin{align}
|F(a)-F(-a)|&\le C_1 a\\
C_2 a &\le \log G(a)-\log G(-a).
\end{align}
We thus find:
\begin{align}
\frac{F(a) G(a)^l - F(-a) G(-a)^l}{a}
&=
G(-a)^l \frac{F(a)-F(-a)}{a}
+
F(a) G(-a)^l \frac{e^{l (\log G(a)-\log G(-a))}-1}{a}\\
&\ge
G(-a)^l
\left(
\frac{F(a)-F(-a)}{a}
+
F(a) \frac{\log G(a)-\log G(-a)}{a}
\right)\\
&\ge
G(-a)^l (F(a) C_2 l-C_1),
\end{align}
which is clearly positive for sufficiently large $l$.
\end{proof}

\section{The main theorem}\label{sec:main}

Given a power series (or left-finite Laurent series) $p$, we define
the valuation $\nu(p)$ of $p$ to be the exponent of the first nonzero
term in $p$.

The general scenario we consider is as follows.  We are given power series
$f(t)$, $g(t)$, and $h(t)$, and asked to prove a statement of the following
form.  If
\[
A(t) = h(t) \sum_{0\le i\le m} c_i f(t)^{m-i} g(t)^i
\]
is such that $A\succeq 0$, then $\nu(A-1)\lesssim \delta m$ as
$m\to\infty$.  (By this, we mean that if $\nu_m$ is the maximum possible
value for given $m$, then $\limsup_{m\to\infty} \nu_m/m \le \delta$.)

We make the normalizing assumptions $f(0)=1$, $h(0)=1$, $g(0)=0$,
$g'(0)=1$.  (We could also assume $f'(0)=0$, but this is somewhat unnatural
in the cases of interest.)

In this context, we recall the following variant of the
B\"urmann-Lagrange theorem:

\begin{lem}
Let $\phi(t)$ and $\psi(t)$ be formal power series, where $\psi(0)=0$,
$\psi'(0)\ne 0$.  Then
\[
\phi(t) = \sum_{0\le i} \kappa_i \psi(t)^i,
\]
where
\[
\kappa_i = [t^0] \phi(t) \frac{t\psi'(t)}{\psi(t)} \psi(t)^{-i}.
\]
\end{lem}

The advantage of this formulation for our purposes is that the dependence
on $j$ below is encoded in a single power series.

\begin{cor}
Let the coefficients $a_j$ and $c_i$ be related by the formal power
series identity
\[
\sum_{j\ge 0} a_j t^j = h(t) \sum_{i\ge 0} c_i f(t)^{m-i} g(t)^i,
\]
with $f$, $g$, $h$ as above.  Then
\[
c_i = \sum_j \alpha_{ij} a_j,
\]
where
\[
\alpha_{ij}
=
[t^0]
t^j h(t)^{-1} f(t)^{i-m} g(t)^{-i}
\left[ \frac{tg'(t)}{g(t)} - \frac{t f'(t)}{f(t)} \right].
\]
\end{cor}

\begin{proof}
We write
\[
t^j = h(t) \sum_{i\ge 0} \alpha_{ij} f(t)^{m-i} g(t)^i,
\]
and thus
\[
h(t)^{-1} f(t)^{-m} t^j = \sum_{i\ge 0} \alpha_{ij} \left(\frac{g(t)}{f(t)}\right)^i;
\]
the formula follows immediately from the lemma.
\end{proof}

In particular,

\begin{cor}
Let $f$, $g$, $h$ be as above, and suppose
\[
\sum_j a_j t^j = h(t) \sum_{0\le i\le m} c_i f(t)^{m-i} g(t)^i.
\]
Then for all $i>m$,
\[
\sum_j \alpha_{ij} a_j = 0.
\]
\end{cor}

We will thus need an asymptotic analysis of the coefficients $\alpha_{ij}$.
In fact, until Section \ref{sec:hermite} we will need only the cases
$i=m+1$ and $i=m+2$.

To make the asymptotic analysis tractable, we need some additional
assumptions (summarized in the statement of Theorem \ref{thm:main1} below).
First, we assume $f\succeq 0$, $f\not\equiv 1$; in all of our applications,
$f$ is derived from a weight enumerator or theta series, so this condition
is automatic.  Our second condition, that $1/g\succeq 0$, is less
automatic, but is easily verified in all cases of interest.

Define a function $Lg(t) = t g'(t)/g(t)$.  By Lemma \ref{lem:logconvex}
applied to the Laurent series $1/g$, we find $Lg'(t)<0$ for positive $t$
within the radius of convergence of $1/g$.  Thus if the equation $Lg(t)=0$
has a positive real root, it must be unique.  We assume the root exists,
and denote it by $t_0$.  Note that if $1/g$ has radius of convergence $r$,
and $\lim_{t\to r^-} g(t)=0$, then $1/g(t)$ is eventually increasing, and
thus $Lg(t)$ is eventually negative; since $Lg(0)=1>0$, this implies that
$Lg(t)=0$ has a positive real root.

We finally assume that $t_0$ is within the open disc of convergence of
$f(t)$, $f(t)/h(t)$, and $f'(t)/h(t)$, and that $1/h(t)>0$ on $[0,t_0]$.

\begin{lem}
Let $f$ and $g$ satisfy the above conditions.  Then $f(t)/g(t)$ has
a unique local minimum on $(0,t_0)$.
\end{lem}

\begin{proof}
Since $f/g\succeq 0$, we conclude that $f(e^s)/g(e^s)$ is strictly convex
on $(-\infty,\log t_0)$.  Now, $\lim_{t\to 0} f(t)/g(t)=+\infty$, and thus
$f(e^s)/g(e^s)$ is decreasing in a neighborhood of $s=-\infty$.  It thus
remains only to show that $f(t)/g(t)$ is increasing near $t_0$.

We compute:
\[
t\frac{d}{dt} \frac{f(t)}{g(t)}
=
g(t)^{-1} (f'(t) - Lg(t) f(t)).
\]
At $t_0$, this becomes $f'(t)/g(t)>0$.
\end{proof}

Denote this local minimum by $t_1$.

\begin{lem}
Let $f$, $g$ satisfy the above assumptions.  Then there exists
a unique point $t'_0$ in $(0,t_0)$ such that
\[
\frac{f(t'_0)}{g(t'_0)} = \frac{f(t_0)}{g(t_0)};
\]
moreover, $t'_0<t_1$.
\end{lem}

\begin{proof}
Uniqueness follows from strict convexity, so it suffices to show existence.
But $f(t)/g(t)$ is continuous and decreasing on $(0,t_1]$, and converges to
$\infty$ at $0$, so attains every value greater than $f(t_1)/g(t_1)$, in
particular $f(t_0)/g(t_0)$.
\end{proof}

The basic idea behind the proof of Theorem \ref{thm:main1} below is that,
to first order, the relation corresponding to
\[
c_{m+2} - \frac{f(t_0)}{g(t_0)} c_{m+1}
\]
is 0 at $j=0$ and is positive for $j>Lg(t'_0)m$.  Moreover, if we
perturb the relation by subtracting a small multiple of $c_{m+1}$, the
relation becomes positive at $j=0$, at the cost of slightly reducing
the range of positivity.  I.e., we can use Lemma \ref{lem:mainlem}
to obtain a contradiction for all relaxations of the desired bound
when $m$ is sufficiently large.

\begin{thm}\label{thm:main1}
Let $f$, $g$, and $h$ be convergent real power series satisfying the following
hypotheses:
\begin{enumerate}
\item $f,h=1+O(t)$, $g=t+O(t^2)$.
\item $f$, $1/g\succeq 0$.
\item $Lg(t):=t g'(t)/g(t)$ has a positive real zero; let the smallest such
zero be denoted $t_0$.
\item $f$, $1/g$, $f/h$, and $f'/h$ have radius of convergence $>t_0$.
\item $1/h$ is positive on $[0,t_0]$.
\end{enumerate}
For each integer $m\ge 0$, let $d_m$ be the maximum of $\nu(A-1)$ where $A$
ranges over power series $\succeq 0$ of the form
\[
A(t)=h(t) \sum_{0\le i\le m} c_i f(t)^{m-i} g(t)^i.
\]
Then $\limsup_{m\to \infty} d_m/m \le Lg(t'_0)$, where $0<t'_0<t_0$ is
such that $f(t'_0)/g(t'_0)=f(t_0)/g(t_0)$.
\end{thm}

\begin{proof}
Choose $0<t_2<t'_0$; we will show that for all sufficiently large $m$,
$d_m/m\le Lg(t_2)$.  As this will hold for all choices of $t_2$, the theorem
will follow from the continuity of $Lg$.  We also choose $0<t_3<t_2$.

We consider the linear combination $c_{m+2}-f(t_3)/g(t_3)c_{m+1}$, as a
vanishing linear combination of the coefficients of $A(t)$.  In particular,
if we let $\alpha_j(t_3)$ denote the coefficient of $[t^j]A(t)$ in this linear
combination, it will suffice to show that for all sufficiently large $m$,
$\alpha_j(t_3)>0$ for $j=0$ and $j\ge Lg(t_2)m$.  Indeed, if a choice of $A(t)$
existed with $d_m\ge Lg(t_2)m$, this would give a positive linear combination
of nonnegative quantities (at least one of which is positive), equal to $0$,
a contradiction.

We compute:
\begin{align}
\alpha_j(t_3)&=\alpha_{(m+2)j}-\frac{f(t_3)}{g(t_3)}\alpha_{(m+1)j}\\
&=
[t^0]
t^j
\left(\frac{f(t)}{g(t)h(t)}\right)
\left(\frac{f(t)}{g(t)}-\frac{f(t_3)}{g(t_3)}\right)
\left(\frac{tg'(t)}{g(t)}-\frac{tf'(t)}{f(t)}\right)
g(t)^{-m}\\
&=
-
[t^{m+2-j}]
t^3
\left(\frac{1}{h(t)}\right)
\left(\frac{f(t)}{g(t)}-\frac{f(t_3)}{g(t_3)}\right)
\frac{d}{dt} \left(\frac{f(t)}{g(t)}\right)
\tilde{g}(t)^{-m},
\end{align}
where $\tilde{g}(t):=t^{-1} g(t)$.  We thus need to show that for all
sufficiently large $m$,
\[
[t^k] F(t) G(t)^m > 0
\]
for $k=m+2$ and $0\le k\le (1-Lg(t_2))m$, where
\begin{align}
G(t) &= \tilde{g}(t)^{-1},\\
F(t) &= t^3 h(t)^{-1}
\left(\frac{f(t_3)}{g(t_3)}-\frac{f(t)}{g(t)}\right)
\frac{d}{dt} \left(\frac{f(t)}{g(t)}\right).
\end{align}

We note the following properties:
\begin{enumerate}
\item $F(t)$ and $G(t)$ are power series with radius of convergence $>t_0$.
\item $G(t)\succeq 0$.
\item $F(t_0)>0$.  Indeed, the first two factors are clearly positive,
while the second factor is positive since
\[
\frac{f(t_0)}{g(t_0)}=\frac{f(t'_0)}{g(t'_0)}<\frac{f(t_3)}{g(t_3)}.
\]
Similarly, the third factor is positive, since $t_0>t_1$.
\item $F(t)$ is positive on $[0,t_3)$.  This time, the second and third
factors are negative; we note the limit $F(0)=1$.
\end{enumerate}

     In other words, the hypotheses of Lemma \ref{lem:mainlem} apply, taking
$I = [0,t_2]\cup [t_0,t_0+\epsilon]$ for appropriate $\epsilon>0$.
It follows that for all sufficiently large $m$,
\[
[t^j] F(t) G(t)^m > 0
\]
when $j/m \in [0,1-Lg(t_2)] \cup [1,1-Lg(t_0+\epsilon)]$.  For sufficiently
large $m$, $(m+2)/m\in [1,1-Lg(t_0+\epsilon)]$, so we are done.
\end{proof}

Example.  Let $A(x,y)$ be the weight enumerator of a doubly-even binary
self-dual code of length $n$.  Then by Gleason's theorem, we have:
\[
A(x,y) = \sum_{0\le i\le [n/24]} c_i
(x^8+14x^4y^4+y^8)^{n/8-3i}(x^4y^4(x^4-y^4)^4)^i.
\]
Defining a power series $A(t)$ by $A(t) = A(1,t^{1/4})$, we obtain:
\begin{align}
A(t) &= \sum_{0\le i\le [n/24]} c_i (1+14t+t^2)^{n/8-3i}(t(1-t)^4)^i\\
     &= (1+14t+t^2)^{n/8-3[n/24]}
       \sum_{0\le i\le [n/24]} c_i ((1+14t+t^2)^3)^{[n/24]-i}(t(1-t)^4)^i.
\end{align}
We will apply the main theorem, with
\begin{align}
f(t) &= (1+14t+t^2)^3\\
g(t) &= t(1-t)^4\\
h(t) &= (1+14t+t^2)^{n/8-3[n/24]}.
\end{align}
We clearly have $f,1/g,h\succeq 0$, and each of $f$, $1/g$, $f/h$, $f'/h$
has radius of convergence at least $1$; indeed, except for $1/g$, these are
polynomials.  Since $\lim_{t\to 1} g(t) = 0$, the hypotheses of Theorem
\ref{thm:main1} are satisfied.

We easily compute $t_0=1/5$; we then find that $t'_0$ is the unique
solution in $(0,1/5)$ of the quartic equation
\[
t^4 - 644 t^3 + 6 t^2 - 644 t + 1.
\]
We thus obtain the bound
\[
\limsup_{m\to\infty} \frac{d_m}{m} \le \frac{1-5 t'_0}{1-t'_0},
\]
or, since $\nu(A(t)-1) = d(C)/4$ and $m = [n/24]$,
\[
\limsup_{n\to\infty} \frac{d(C)}{n(C)} \le \frac{1-5 t'_0}{6(1-t'_0)}.
\]
In fact (as we will explain below), $t'_0$ has the simple closed form
\[
t'_0 = \left(\frac{1-5^{-1/4}}{1+5^{-1/4}}\right)^4;
\]
when substituted in, this simplifies (again explained below) to give
\[
\limsup_{n\to\infty} d_n/n \le \frac{1}{2} (1-5^{-1/4}).
\]

\section{Applications I}\label{sec:apps1}

We generalize the previous example as follows (compare the Gleason-Pierce
theorem \cite{SloaneNJA:1979}):

\begin{thm}\label{thm:codes2}
Let $q$ and $c$ be chosen with $q>1$, such that either $c=1$ or $(q,c)\in
\{(2,2),(2,4),(3,3),(4,2)\}$.  Let $C_i$ be a sequence of formally
self-dual codes over an alphabet of size $q$ with all Hamming weights
divisible by $c$; suppose furthermore that as $i\to\infty$,
$n(C_i)\to\infty$.  Then
\[
\limsup_{i\to\infty} \frac{d(C_i)}{n(C_i)}
\le
\frac{q-1}{q}(1-(c+1)^{-1/c}).
\]
\end{thm}

\begin{proof}
Given such a code $C$, of length $n$, let $A_C(x,y)$ be its weight
enumerator, and define a power series $A(t)=A_C(1,t^{-1/c})$.  Then from
the various Gleason theorems \cite[Section 7]{me:N}, we conclude
\begin{align}
c=1: A(t) &= (1+(\sqrt{q}-1)t)^{n\bmod 2} \sum_{0\le i\le [n/2]}
c_i ((1+(\sqrt{q}-1)t)^2)^{[n/2]-i} (t(1-t))^i\\
(q,c)=(2,2): A(t) &= (1+t)^{(n/2)\bmod 4}
\sum_{0\le i\le [n/8]} c_i ((1+t)^4)^{[n/8]-i} (t(1-t)^2)^i\\
(q,c)=(2,4): A(t) &= (1+14t+t^2)^{(n/8)\bmod 3}
\sum_{0\le i\le [n/24]} c_i ((1+14t+t^2)^3)^{[n/24]-i} (t(1-t)^4)^i\\
(q,c)=(3,3): A(t) &= (1+8t)^{(n/4)\bmod 3}
\sum_{0\le i\le [n/12]} c_i ((1+8t)^3)^{[n/12]-i} (t(1-t)^3)^i\\
(q,c)=(4,2): A(t) &= (1+3t)^{[n/2]\bmod 3}
\sum_{0\le i\le [n/6]} c_i ((1+3t)^3)^{[n/6]-i} (t(1-t)^2)^i,
\end{align}
for appropriate coefficients $c_i$.
In particular, we are in the scenario of Theorem \ref{thm:main1},
with $g(t) = t(1-t)^c$.  In each case, $f$, $f/h$ and $f'/h$
are all clearly polynomials; since $1/g\succeq 0$, has radius of
convergence $1$, and $g(1)=0$, the hypotheses of the theorem are
satisfied; it remains to compute $Lg(t'_0)$.

Since
\[
Lg(t) = \frac{1-(1+c)t}{1-t},
\]
we find $t_0=1/(1+c)$.  To compute $t'_0$ from $t_0$, we proceed as follows.

Define new series $F(t)=f(t^c)$ and $G(t)=g(t^c)$.
From the MacWilliams identity and the fact that $f$ and $g$
are linear combinations of power series coming from weight enumerators,
we find that there exists an integer $n_0$ such that:
\begin{align}
F\left(\frac{1-t}{1+(q-1)t}\right)
&=
q^{n_0/2} (1+(q-1)t)^{-n_0} F(t)\\
G\left(\frac{1-t}{1+(q-1)t}\right)
&=
q^{n_0/2} (1+(q-1)t)^{-n_0} G(t).
\end{align}
Dividing these equations, we find
\[
\frac{F}{G}\left(\frac{1-t}{1+(q-1)t}\right)
=
\frac{F}{G}\left(t\right).
\]
In terms of $f$ and $g$, this becomes:
\[
\frac{f}{g}\left(\left(\frac{1-t^{1/c}}{1+(q-1)t^{1/c}}\right)^c\right)
=
\frac{f}{g}(t).
\]
We thus conclude that
\[
t'_0 = \left(\frac{1-(c+1)^{-1/c}}{1+(q-1)(c+1)^{-1/c}}\right)^c,
\]
since we readily verify $0<t'_0<t_0$.

Similarly, to compute $Lg(t'_0)$, we differentiate the functional equation
for $G$ at $t_0^{1/c}$, obtaining:
\[
Lg(t'_0)
=
\frac{(q-1)n_0}{cq} (1-t_0^{1/c})
+
Lg(t_0) \left(\frac{(q-1)t_0^{1/c} - (q-2) - t_0^{-1/c}}{q}\right)
=
\frac{(q-1)n_0}{cq} (1-t_0^{1/c}),
\]
since by definition $Lg(t_0)=0$.  Since $m = [n/n_0]$ and
$\nu(A(t)-1) = d(C)/c$, we mutiply this by $c/n_0$ to obtain
the desired bound.
\end{proof}

\begin{rem}
In particular:
\begin{enumerate}
\item For doubly-even self-dual binary codes, $d/n\lesssim (1-5^{-1/4})/2 =
.1656298476$.
\item For self-dual ternary codes, $d/n\lesssim (2-2^{1/3})/3 =
.2466929834$.
\item For even self-dual additive codes over $GF(4)$, $d/n\lesssim
(3-3^{1/2})/4 = .3169872982$.
\item For singly-even, formally self-dual binary codes, $d/n\lesssim
(1-3^{-1/2})/2 = .2113248655$.
\item For formally self-dual codes over $GF(q)$, $d/n\lesssim 1/2-1/2q$.
\end{enumerate}
Strictly speaking only the case $(q,c)=(3,3)$ is new; the case
$(q,c)=(2,4)$ was shown (via a rather different proof) in
\cite{KrasikovI/LitsynS:2000}, while in the remaining cases, the
bound obtained is worse than the JPL
\cite{McElieceRJ/RodemichER/RumseyH/WelchL:1977} or Aaltonen
\cite{AaltonenMJ:1979}, \cite{AaltonenM:1990} bound, as appropriate.
But for singly-even self-dual binary codes, see Theorem \ref{thm:codes1}
below.
\end{rem}

For even modular lattices (see \cite{QuebbemannH-G:1995},
\cite{QuebbemannH-G:1997}, \cite{me:M}; note that a unimodular lattice is
1-modular), we proceed similarly.  We recall Dedekind's $\eta$ function
\[
\eta(z) = e^{\pi i z/12}\prod_i (1-e^{2\pi i m z})
\]
and the Eisenstein series
\begin{align}
E_2(z)
&= \frac{1}{2\pi i} \frac{d}{dz} \log\eta(z)\\
&=\frac{1}{24}-\sum_{1\le m} (\sum_{k|m} k) e^{2\pi i m z}.
\end{align}

\begin{thm}\label{thm:latts2}
Let $N$ be one of the integers $\{1,2,3,5,6,7,11,14,15,23\}$, and define
\[
E^{(N)}_2(z) = \sum_{m|N} m E_2(mz).
\]
Then for any sequence $\Lambda_i$ of even, strongly $N$-modular lattices of
dimension tending to $\infty$,
\[
\limsup_{i\to\infty}
\frac{\mu(\Lambda_i)}{\dim(\Lambda_i)}
\le
\frac{N z_0}{2\pi i},
\]
where $z_0$ is the unique zero of $E^{(N)}_2$ on the positive imaginary
axis.
\end{thm}

\begin{proof}
Let $\Lambda$ be an even, strongly $N$-modular lattice, with theta series
$\Theta_\Lambda$.  Then
(\cite{QuebbemannH-G:1995},\cite{QuebbemannH-G:1997})
$\Theta_\Lambda(e^{\pi i z})$ can be written as a weighted-homogeneous
polynomial in $\Theta_{\Lambda_0}(e^{\pi i z})$ and $g(e^{2\pi i z}) =
(\prod_{m|N}
\eta(mz))^{24/(\sum_{m|N} m)}$, where $\Lambda_0$ is the lowest-dimensional
even $N$-modular lattice.

Clearly $\Theta_{\Lambda_0}(t^{1/2})\succeq 0$; using the product
formula for $\eta$, we also conclude that $1/g(t)\succeq 0$.
Thus Theorem \ref{thm:main1} applies, and it remains only to compute
$Lg(t'_0)$.

We first compute, with $t=e^{2\pi i z}$,
\begin{align}
Lg(z) &= \frac{1}{2\pi i} \frac{d}{dz} \log\prod_{m|N}
\eta(mz)^{24/(\sum_{m|N} m)}\\
&=
\frac{24}{\sum_{m|N} m} E^{(m)}_2(z),
\end{align}
and thus $t_0 = e^{2\pi i z_0}$.

We have the transformation laws:
\begin{align}
f(-1/Nz) &= N^{n_0/4} (z/i)^{n_0/2} f(z)\\
g(-1/Nz) &= N^{n_0/4} (z/i)^{n_0/2} g(z),
\end{align}
where
\[
n_0 = \frac{24\delta(N)}{\sum_{m|N} m}.
\]
We thus conclude that $t'_0 = e^{-2\pi i/N z_0}$, and that
\[
Lg(t'_0)
=
\frac{z_0 N n_0}{4\pi i}
+
N z_0^2 Lg(t_0)
=
\frac{z_0 N n_0}{4\pi i}.
\]
Multiplying by $2$ (since the lattices are even) and dividing by $n_0$
(since $m = [n/n_0]$) gives the required bound.
\end{proof}

\begin{rem}
Numerically, we have:
\begin{align}\notag
N= 1:& \frac{\mu(\Lambda)}{\dim(\Lambda)} \lesssim .0833210664&
N= 2:& \frac{\mu(\Lambda)}{\dim(\Lambda)} \lesssim .1246710056\\\notag
N= 3:& \frac{\mu(\Lambda)}{\dim(\Lambda)} \lesssim .1643714543&
N= 5:& \frac{\mu(\Lambda)}{\dim(\Lambda)} \lesssim .2351529896\\\notag
N= 6:& \frac{\mu(\Lambda)}{\dim(\Lambda)} \lesssim .2414115212&
N= 7:& \frac{\mu(\Lambda)}{\dim(\Lambda)} \lesssim .2957105217\\\notag
N=11:& \frac{\mu(\Lambda)}{\dim(\Lambda)} \lesssim .3973198712&
N=14:& \frac{\mu(\Lambda)}{\dim(\Lambda)} \lesssim .4266498017\\\notag
N=15:& \frac{\mu(\Lambda)}{\dim(\Lambda)} \lesssim .3206725342&
N=23:& \frac{\mu(\Lambda)}{\dim(\Lambda)} \lesssim .6262824896
\end{align}
Again, aside from $N=1$, $N=2$, $N=3$, the obtained bound is worse than
that implied by the Kabatiansky-Levenshtein bound on sphere packings
\cite{KabatianskyGA/LevenshteinVI:1978}.
\end{rem}

We finally consider self-dual codes over $\Z_4$.  As in \cite{me:CC},
bounding the Hamming or Lee distance reduces to a consideration of the dual
distance of doubly-even binary codes.  At length a multiple of 8, the bound
on self-dual doubly-even codes applies; for other lengths, we shorten the
code up to 7 times, without affecting the asymptotic bound.  We thus
obtain:

\begin{thm}\label{thm:Z4HL}
Let $C_i$ be a sequence of self-dual codes over $\Z_4$, with length
tending to $\infty$.  Then
\begin{align}
\limsup_{i\to\infty} \frac{d_H(C_i)}{n(C_i)}
&\le
\frac{1-5^{-1/4}}{2}\\
\limsup_{i\to\infty} \frac{d_L(C_i)}{n(C_i)}
&\le
1-5^{-1/4},
\end{align}
where $d_H(C)$ and $d_L(C)$ are the minimum Hamming and Lee weights of
$C$, respectively.
\end{thm}

For the Euclidean distance, the situation is more complicated.
We use the following lemma:

\begin{lem}
Let $C$ be a Type II self-dual code over $\Z_4$ (all Euclidean norms
divisible by $8$), and let $A(x,y,z)$ be its symmetrized weight enumerator.
Define a power series
\[
A(t):=A(1,2^{-1/4} t^{1/8}(1+t)^{1/4},t^{1/2}).
\]
Then $A\succeq 0$, and
\[
\nu(A(t)-1) = d_E(C)/8.
\]
Moreover, for appropriate coefficients $c_i$,
\[
A(t) =
\sum_{0\le i\le [n/24]}
c_i
(1+60t+134 t^2+60 t^3+t^4)^{(n/8)-3i}
(t (1-t)^6 (1-t^2)^2)^i.
\]
\end{lem}

\begin{proof}
Consider a monomial $x^a y^b z^c$ of $A(x,y,z)$; note that
$b+4c$ must be a multiple of 8, and in particular $b$ is a multiple
nof 4.  Under the specified substitution, this yields
\[
2^{-b/4} (1+t)^{b/4} t^{(b+4c)/8},
\]
a polynomial with nonnegative coefficients and with valuation $(b+4c)/8$;
the first two claims are immediate.

Now, the Gleason theorem for Type II codes over $\Z_4$ states that
$A(x,y,z)$ is a (weighted homogeneous) polynomial in
\begin{align}
\theta_8 &= \frac{(x+z)^8+(x-z)^8}{2} + 128 y^8\\
\theta_{16} &= (x^2z^2(x^2+z^2)^2-4 y^8)((x^4+6x^2z^2+z^4)^2-64y^8)\\
\theta_{24} &= y^8 (x^2-z^2)^8\\
h_8 &= (xz(x^2+z^2)-2y^4)^2.
\end{align}
Under the substituion, we have:
\begin{align}
\theta_8&\mapsto 1+60t+134 t^2+60 t^3+t^4\\
\theta_{16}&\mapsto 0\\
\theta_{24}&\mapsto \frac{1}{4} t (1-t)^6 (1-t^2)^2\\
h_8 &\mapsto 0
\end{align}
The remaining claim follows.
\end{proof}

\begin{rem}
The above substitution is inspired by the proofs used in
\cite{BonnecazeA/SoleP/BachocC/MourrainB:1997} and
\cite{me:M}, which involve lifting the code to a lattice and analyzing
the resulting theta series.  The resulting substitution takes the
polynomial
\[
2y^4 - xz(x^2+z^2)
\]
to 0; solving gives the above substitution.
\end{rem}

We apply Theorem \ref{thm:main1} with $f(t) = (1+60t+134 t^2+60
t^3+t^4)^3$, $g(t) = t(1-t)^6(1-t^2)^2$, to obtain:

\begin{thm}\label{thm:Z4E2}
Let $C_i$ be a sequence of Type II self-dual codes over $\Z_4$ with
length tending to $\infty$.  Then
\[
\limsup_{i\to\infty} \frac{d_E(C_i)}{n(c_i)}
\le .3332625492 = 2-2x,
\]
where $x$ is the positive real root of the polynomial
\[
11 x^{16}+2112 x^{14}-8525 x^{12}+15048 x^{10}-15218 x^8+9552 x^6-3718 x^4+828 x^2-81.
\]
\end{thm}

\begin{rems}
The field $\Q[x]$ is that generated by
\[
2^{-1/4} t_0^{1/8}(1+t_0)^{1/4}\text{ and }t_0^{1/2},
\]
where $t_0 = \frac{2\sqrt{5}-3}{11}$.
\end{rems}

\begin{rems}
The Type II hypothesis is removed in Theorem \ref{thm:Z4E1} below.
\end{rems}

\section{Applications II}\label{sec:apps2}

We now consider applications to which Theorem \ref{thm:main1} does not
directly apply, but for which the same basic idea can be used.

We first extend the bound for doubly-even self-dual binary codes to
general self-dual binary codes.  The main idea is that, using Lemma
\ref{lem:shadow}, we can reduce positivity of coefficients of the
form
\[
[t^j] F_1(t) F_2(t)^l G(t^2)^m
\]
to positivity of $F_1(t)$, for all sufficiently large $l$ and all
sufficiently large $m$; this because
\begin{align}
[t^{2j}] F_1(t) F_2(t)^l G(t^2)^m
&=
[t^j] \left(F_1(\sqrt{t})F_2(\sqrt{t})^l +
F_1(-\sqrt{t})F_2(-\sqrt{t})^l\right) G(t)^m\\
[t^{2j+1}] F_1(t) F_2(t)^l G(t^2)^m
&=
[t^j] \left(\frac{F_1(\sqrt{t})F_2(\sqrt{t})^l -
F_1(-\sqrt{t})F_2(-\sqrt{t})^l}{\sqrt{t}}\right) G(t)^m.
\end{align}

\begin{thm}\label{thm:codes1}
Let $C_i$ be a sequence of self-dual binary codes of length tending to
$\infty$.  Then
\[
\limsup_{i\to\infty} \frac{d(C_i)}{n(C_i)}
\le
\frac{1-5^{-1/4}}{2}.
\]
\end{thm}

\begin{proof}
We recall that a self-dual binary code $C$ has associated to it two
enumerators: its weight enumerator $A(x,y)$ and its ``shadow'' enumerator
$S(x,y)$ \cite{ConwayJH/SloaneNJA:1990}, \cite{me:N}.\footnote{Note that
this includes the
doubly-even case, for which $S(x,y)=A(x,y)$.}  For appropriate coefficients
$c_i$, we have:
\begin{align}
A(x,y) &=
(x^2+y^2)^{(n\bmod 8)/2}
  \sum_{0\le i\le [n/8]} c_i ((x^2+y^2)^4)^{[n/8]-i} (x^2y^2(x^2-y^2)^2)^i\\
S(x,y) &=
(2xy)^{(n\bmod 8)/2}
  \sum_{0\le i\le [n/8]} c_i ((2xy)^4)^{[n/8]-i} (-(x^4-y^4)^2/4)^i
\end{align}
Here $c_i$ can be characterized (up to a multiplicative constant) as the
unique linear combination of the first $i+1$ nontrivial coefficients of
$A(x,y)$ that is also a linear combination of the first $[n/8]-i$
nontrivial coefficients of $S(x,y)$.

We define power series $A(t)$ and $S(t)$ by
\begin{align}
A(t) &= A(1,t^{1/2})\\
     &= (1+t)^{(n\bmod 8)/2} \sum_{0\le i\le [n/8]} c_i ((1+t)^4)^{[n/8]-i}
                                                        (t(1-t)^2)^i\\
S(t) &= t^{-(n\bmod 8)/8} S(1,t^{1/4})\\
     &= 2^{(n\bmod 8)/2} 
        \sum_{0\le i\le [n/8]} c_i (16t)^{[n/8]-i} (-(1-t)^2/4)^i.
\end{align}
If we let $\alpha_{ij}$ denote the coefficient of $[t^j] A(t)$ in
$c_i$, and let $\beta_{ij}$ denote the coefficient of $[t^j] S(t)$ in
$c_i$, B\"urmann-Lagrange tells us:
\begin{align}
\alpha_{ij} &= [t^{i-j}] (1-6t+t^2) (1+t)^{4i-n/2-1} (1-t)^{-2i-1}\\
            &= [t^{i-j}] (1-6t+t^2) (1+t)^{6i-n/2} (1-t^2)^{-2i-1}\\
\beta_{ij}  &= (-1)^i 2^{6i-n/2} [t^{[n/8]-i-j}] (1+t) (1-t)^{-2i-1}.
\end{align}

Set $t_0=5^{-1/2}$, $t'_0 = \left(\frac{1-5^{-1/4}}{1+5^{-1/4}}\right)^2$,
and let $m$ be an even integer of the form $n/12+l/6$.  Now, consider the
coefficient
\[
c_{m+1}-\frac{(1+t_0)^4}{t_0(1-t_0)^2} c_m.
\]
This is manifestly a negative linear combination of $[t^j] S(t)$ for $0\le
j\le [n/8]-m$ (so the same will be true for small perturbations).
It thus suffices to consider the coefficients of $[t^j] S(t)$.
By the remark above and Lemma \ref{lem:shadow}, this reduces to showing that
\[
(1-6t+t^2) \left( \frac{(1+t)^4}{t(1-t)^2} -
\frac{(1+t_0)^4}{t_0(1-t_0)^2}\right)
\]
is positive for $0\le t<t'_0$, and can be perturbed to be positive at $t_0$
as well; this is clearly the case.
\end{proof}

Similarly, the bound for even self-dual additive codes over $GF(4)$
extends to the general case; the resulting bound is still weaker than the
Aaltonen bound, however.  A similar remark applies to formally self-dual
binary codes, with shadow defined by $S(x,y)=A(y,x)$; there the bound
for singly-even f.s.d. binary codes extends.  In that case, $m\sim n/4$,
instead of the obvious analogue $n/3$; also we must apply Lemma
\ref{lem:shadow} to the coefficients of $S$ as well as to the coefficients
of $A$.

For lattices, we have:

\begin{thm}\label{thm:latts1}
Let $N$ be one of the integers $\{1,2,3,5,6,7,11,14,15,23\}$, and define
\[
E^{(N)}_2(z) = \sum_{m|N} m E_2(mz).
\]
Then for any sequence $\Lambda_i$ of strongly $N$-modular lattices of
dimension tending to $\infty$,
\[
\limsup_{i\to\infty}
\frac{\mu(\Lambda_i)}{\dim(\Lambda_i)}
\le
\frac{N z_0}{2\pi i},
\]
where $z_0$ is the unique zero of $E^{(N)}_2$ on the positive imaginary
axis.
\end{thm}

\begin{proof}
As above; the case $N$ odd is analogous to the self-dual binary code case,
while the case $N$ even is analogous to the formally self-dual binary code
case.  The only respect in which the proof is not straightforward (using
the formulae of \cite{me:M}) is in dealing with the ``other'' genera
(not covered by Theorem 2 of \cite{me:M}).  In each case, direct
summing by a suitable $N$-modular lattice places us into the ``good''
genera, and we can proceed from there; the only effect is to multiply the
power series in question by a theta series, which clearly has no effect on
positivity.
\end{proof}

Our last shadow application is to codes over $\Z_4$:

\begin{thm}\label{thm:Z4E1}
Let $C_i$ be a sequence of self-dual codes over $\Z_4$ with
length tending to $\infty$.  Then
\[
\limsup_{i\to\infty} \frac{d_E(C_i)}{n(c_i)}
\le .3332625492 = 2-2x,
\]
where $x$ is the positive real root of the polynomial
\[
11 x^{16}+2112 x^{14}-8525 x^{12}+15048 x^{10}-15218 x^8+9552 x^6-3718 x^4+828 x^2-81.
\]
\end{thm}

\begin{proof}
We define
\begin{align}
A(t) &= A(1,2^{-1/4} t^{1/4}(1+t^2)^{1/4},t)\\
S(t) &= t^{-(n\bmod 8)/8}(1+t)^{-(n\bmod 4)/4}
        S(1,2^{-1/4} t^{1/8}(1+t)^{1/4},t^{1/2}),
\end{align}
and observe that $A$ and $S$ are polynomials with nonnegative coefficients.
Since
\begin{align}
A(t) &= 
\sum_{0\le i\le [n/8]} c_i ((1+t)^8)^{n/8-i} (t(1+t^2)(1-t)^4)^i\\
S(t) &= 
(1+t)^{[n/4] \bmod 2}
\sum_{0\le i\le [n/8]} (-1)^i (64 t (1+t)^2)^{[n/8]-i} ((1-t)^4/8)^i,
\end{align}
we can proceed as in the case of self-dual binary codes.
\end{proof}

Our final application is to quantum codes.  Formally, a $q$-ary quantum
code corresponds to a self-orthogonal codes $C$ over an alphabet of size
$q^2$; the objective is to bound the minimum weight of the nonlinear code
$C^\perp-C$.  We recall the following (the nonbinary extension of
\cite[Theorem 6]{me:X}):

\begin{lem}
Let $Q$ be a $q$-ary quantum code of length $n$ and dimension $K$.  Then
there exist polynomials $C(x,y)$ and $D(x,y)$, homogeneous of degree $n$,
such that
\begin{align}
C(x,y) &= C(\frac{x+(q^2-1)y}{q},\frac{x-y}{q})\\
D(x,y) &= -D(\frac{x+(q^2-1)y}{q},\frac{x-y}{q}),\\
C(1,0) &= 1,\\
\intertext{and satisfying}
\nu(C(1,t)-D(1,t)) &= d\\
C(1,t) &\succeq \frac{K-1}{2K} (C(1,t)-D(1,t))\label{eq:Cconstr}\\
C(1,t)-D(1,t) &\succeq 0.
\end{align}
\end{lem}

In particular, if $C$ and $D$ satisfy the constraints for a given value
$K>1$, they satisfy them for all smaller $K$, including for $K=1$.
Thus if we replace \eqref{eq:Cconstr} by the condition
\[
C(1,t)\succeq 0,
\]
the resulting bound will apply to all quantum codes.

\begin{thm}\label{thm:quantq}
Let $Q_i$ be a sequence of $q$-ary quantum codes of length tending to
$\infty$.  Then
\[
\limsup_{i\to\infty} \frac{d(Q_i)}{n(Q_i)}
\le
\frac{1}{2}\left(1-\frac{1}{q^2}\right).
\]
\end{thm}

\begin{proof}
We consider the case of odd length $n=2m+1$; the case of even length is
analogous.  Setting $C(t):=C(1,t)$, $D(t)=D(1,t)$, we observe that
\begin{align}
C(t) &=(1+(q-1)t) \sum_{0\le i\le m} c_i ((1+(q-1)t)^2)^{m-i} (t(1-t))^i\\
D(t) &=(1-(q+1)t) \sum_{0\le i\le m} d_i ((1+(q-1)t)^2)^{m-i} (t(1-t))^i,
\end{align}
for suitable coefficients $c_i$ and $d_i$.  Let $\gamma_{ij}$ be the
coefficient of $[t^j] C(t)$ in $c_i$ (extending as usual to $i>m$), and let
$\delta_{ij}$ be the coefficient of $[t^j] D(t)$ in $d_i$.  We find:
\begin{align}
\gamma_{(m+1)j} &= [t^{m+2-j}] (1-(q+1)t) t(1-t) (1-t)^{-m-3}\\ 
\gamma_{(m+2)j} &= [t^{m+2-j}] (1-(q+1)t) (1+(q-1)t)^2 (1-t)^{-m-3}\\ 
\delta_{(m+1)j} &= [t^{m+2-j}] (1+(q-1)t) t(1-t) (1-t)^{-m-3}\\ 
\delta_{(m+2)j} &= [t^{m+2-j}] (1+(q-1)t) (1+(q-1)t)^2 (1-t)^{-m-3}.
\end{align}
If we instead expand $c_i$ and $d_i$ in terms of $[t^j] C(t)$ and
$[t^j] (C(t)-D(t))$, we obtain coefficients:
\begin{align}
\left(\text{coeff. of $[t^j] C(t)$ in $c_i$}\right) &= \gamma_{ij}\\
\left(\text{coeff. of $[t^j] (C(t)-D(t))$ in $c_i$}\right) &= 0\\
\left(\text{coeff. of $[t^j] C(t)$ in $d_i$}\right) &= \delta_{ij}\\
\left(\text{coeff. of $[t^j] (C(t)-D(t))$ in $d_i$}\right) &= -\delta_{ij}
\end{align}
We need a relation that is a nonnegative linear combination of
the coefficients $[t^j] C(t)$, positive at $j=0$, as well as a
nonnegative linear combination of the coefficients $[t^j] C(t)-D(t)$
for $j$ larger than the bound.

Now, consider the relation
\[
(q+1)(c_1 - (q+1)^2 c_0)-(q-1)(d_1-(q+1)^2 d_0).
\]
On $[t^j] C(t)-D(t)$, this has coefficient:
\[
(q-1)[t^{m+2-j}] (1+(q-1)t) (1-2t) (1-(q^2+1)t) (1-t)^{-m-3},
\]
which to first order is positive for $j/m>1-\frac{1}{q^2}$.
Similarly, on $[t^j] C(t)$, this has coefficient:
\[
2 [t^{m+2-j}] (1-2t) (1-(q^2+1)t)^2 (1-t)^{-m-3}
\]
which is positive for $0\le j\le m+2$ except in a neighborhood of
$j=0$ and $j=(1-\frac{1}{q^2})m$.  The construction of a positive
perturbation is straightforward.
\end{proof}

Similarly, taking the shadow constraints \cite{me:W} into account, we obtain:

\begin{thm}\label{thm:quant2}
Let $Q_i$ be a sequence of binary quantum codes of length tending to
$\infty$.  Then
\[
\limsup_{i\to\infty} \frac{d(Q_i)}{n(Q_i)}
\le
\frac{3-\sqrt{3}}{4}.
\]
\end{thm}

Note that although as we have remarked, this is slightly worse than the
Aaltonen bound, this is still a new result; in the quantum case,
the Aaltonen bound is only known for a set of rates bounded above 0
\cite{AshikhminA/LitsynS:1999}.

\section{Extensions}\label{sec:hermite}

As we have remarked, many of our bounds are weaker than the appropriate
``universal'' bounds (JPL \cite{McElieceRJ/RodemichER/RumseyH/WelchL:1977};
Aaltonen \cite{AaltonenMJ:1979}, \cite{AaltonenM:1990};
Kabitiansky-Levenshtein \cite{KabatianskyGA/LevenshteinVI:1978}; see
\cite{LevenshteinVI:1998} for a survey) that hold even for non-self-dual
codes of rate $1/2$ and non-lattice packings of appropriate density.  Since
others of our bounds are quite a bit stronger than the corresponding
universal bounds, this strongly suggests that in no case is either bound
tight for self-dual codes.  The question then becomes that of how to
improve the above bounds.

We restrict our attention to the situation of Section \ref{sec:main};
we will comment on the shadow and quantum cases at the end.

Thus, let $f$, $g$, $h$ satisfy the hypotheses of Theorem \ref{thm:main1}.
As the above bounds resulted from considering the two relations $c_{m+1}$,
$c_{m+2}$ in place of the single relation $c_{m+1}$, the obvious thing to
try is a linear combination of $c_{m+i}$ for $1\le i\le l$.  Somewhat
surprisingly, this does not give rise to any improvement in terms
of $\limsup d/m$.  We can see this as follows.  The coefficient
of $[t^j] A(t)$ in such a linear combination will have the form
\[
-
[t^{m-j}]
\left(\frac{1}{h(t)}\right)
p(\frac{f(t)}{g(t)})
t\frac{d}{dt} \left(\frac{f(t)}{g(t)}\right)
\tilde{g}(t)^{-m},
\]
for some polynomial $p$.  The point, then, is that since
\[
\frac{d}{dt} \left(\frac{f(t)}{g(t)}\right)
\]
has opposite signs at $t=t_0$ and at $t=t'_0$, while
\[
\left(\frac{1}{h(t)}\right)
\text{ and }
p(\frac{f(t)}{g(t)})
\]
both have the same sign at the two points (the same value in the latter
case), the corresponding coefficients will, to first order, also have
opposite sign.  In particular, we will never obtain a bound on $d/m$ better
than $Lg(t'_0)$.  (This tends to explain why Krasikov and Litsyn obtained
the same bound in the doubly-even binary case, despite a rather different
argument, and the unlikelihood that the bound is optimal.)

On the other hand, we do obtain a slight lower-order improvement:

\begin{thm}\label{thm:hermite}
Let $f$, $g$, $h$, $d_m$, $t_0$, $t'_0$ be as in the hypotheses of Theorem
\ref{thm:main1}, and suppose further that
\[
\frac{t \left(\frac{d}{dt} \frac{f(t)}{g(t)}\right)^2}
{-t Lg'(t)}
\]
is smaller at $t_0$ than at $t'_0$.  Then
\[
\limsup_{m\to\infty}
\frac{d_m-Lg(t'_0)m}{\sqrt{m}}
=
-\infty.
\]
\end{thm}

\begin{proof}
For each $k\ge 0$, let $c^{(k)}$ be the relation corresponding to the
polynomial $p(t) = (t-f(t_0)/g(t_0))^k$ above; that is:
\[
c^{(k)}
=
\sum_{0\le i\le k}
(-1)^i
\binom{k}{i}
\left(\frac{f(t_0)}{g(t_0)}\right)^i
c_{m+1+i}.
\]
Thus the coefficient of $[t^j] A(t)$ in $c^{(k)}$ is
\[
\alpha^{(k)}_j:=
-
[t^{m-j}]
\left(\frac{1}{h(t)}\right)
t\frac{d}{dt} \left(\frac{f(t)}{g(t)}\right)
(\frac{f(t)}{g(t)}-\frac{f(t_0)}{g(t_0)})^k
\tilde{g}(t)^{-m}.
\]
Also, define coefficients $a^{(k)}_j$ by
\[
a^{(k)}_j := \frac{\alpha^{(k)}_j}{[t^{m-j}] \tilde{g}(t)^{-m}}.
\]

We will consider relations of the form
\[
\sum_{0\le k\le l} b_k m^{k/2} c^{(k)}
\]
with $b_l>0$.  We first claim that for any such relation, the coefficients
of $[t^j] A(t)$ are nonnegative when $j\ge Lg(t'_0) n + n^{k/(k+1)}$ and $n$ is
sufficiently large.  Indeed, in a neighborhood of $j/n=1$, this follows
from the estimate of Lemma \ref{lem:mainlem}; the terms for $k<l$ are
$o(m^{l/2})$, while the term for $k=l$ is $\Omega(m^{l/2})$ and positive.
In the other region, we use the $O(n^{-1})$ estimate of Lemma
\ref{lem:asymptseries}, in which the main term is positive and of order
$\Omega(n^{-k/(k+1)})>O(n^{-1})$.

In the remaining neighborhood of $t'_0$, we use Corollary
\ref{cor:asympt_hermite}.  We thus have:
\[
a^{(k)}_{Lg(t'_0)m+x\sqrt{m}}
=
C'_1
(C'_2)^{k+1}
(C'_3 m)^{-k/2}
h_k(x/\sqrt{C'_3})
+
o(m^{-k/2}\max(|x|^k,1)),
\]
where
\begin{align}
C'_1
&=
h(t'_0)^{-1}>0\\
C'_2
&=
-t'_0\frac{d}{dt} \left(\frac{f(t)}{g(t)}\right)_{t=t'_0}>0\\
C'_3
&=
-t'_0 Lg'(t'_0)>0.
\end{align}
In particular, if $x$ is bounded above the largest zero of
\[
\sum_{0\le k\le l} b_k C'_1 (C'_2)^{k+1} (C'_3)^{-k/2} h_k(x/\sqrt{C'_3}),
\]
then we have positivity for $j\ge Lg(t'_0)m+x\sqrt{m}$;
when $j\le Lg(t'_0)m + m^{-l/(l+1)}$, the error term is uniformly of
smaller order than the main term.

Similarly,
\[
a^{(k)}_{0}
=
C_1
(C_2)^{k+1}
(C_3 m)^{-k/2}
h_k(0)
+
o(m^{-k/2}),
\]
with
\begin{align}
C_1
&=
h(t_0)^{-1}>0\\
C_2
&=
-t_0\frac{d}{dt} \left(\frac{f(t)}{g(t)}\right)_{t=t_0}<0\\
C_3
&=
-t_0 Lg'(t_0)>0.
\end{align}

Thus given any choice of $b_k$ such that
\[
q:=\sum_{0\le k\le l} b_k C_1 (C_2)^{k+1} (C_3)^{-k/2} h_k(0) > 0,
\]
if $x(b)$ is the largest zero of
\[
p(x):=\sum_{0\le k\le l} b_k C'_1 (C'_2)^{k+1} (C'_3)^{-k/2}
h_k(x/\sqrt{C'_3}),
\]
we have the asymptotic bound
\[
\limsup_{m\to\infty} \frac{d_m-Lg(t'_0)m}{m^{1/2}} \le x(b).
\]

To construct a good relation, we will need some further properties of
the Hermite polynomials, all classical results:

\begin{lem}
The polynomials $h_k(x)$ are the unique monic polynomials such that
\[
\frac{1}{\sqrt{2\pi}} \int_{-\infty}^\infty
h_j(x) h_k(x) e^{-x^2/2} dx
=
\delta_{jk} k!.
\]
Furthermore, we have the three-variable generating function:
\[
\sum_{k\ge 0} h_k(x) h_k(y) \frac{t^k}{k!}
=
\frac{1}{\sqrt{1-t^2}}
\exp\left(
-\frac{t^2 x^2/2-txy+t^2y^2/2}{1-t^2}\right),
\]
convergent for $|t|<1$.  Finally (Christoffel-Darboux):
\[
\frac{h_k(x)h_{k-1}(y)-h_{k-1}(x)h_k(y)}{x-y}
=
\sum_{0\le j<k} \frac{h_j(x) h_j(y)}{j!}.
\]
\end{lem}

Now, given $p(x)$, we can compute $b_k$ using orthogonality; we find:
\[
b_k
=
\frac{1}{C'_1 C'_2}
\left(\frac{C'_2}{\sqrt{C'_3}}\right)^{-k}
\frac{1}{\sqrt{2\pi}}
\frac{1}{k!}
\int p(x\sqrt{C'_3}) h_k(x) e^{-x^2/2} dx;
\]
using the same formula to define $b_k$ for $k>l$ gives $b_k=0$.
But then $q$ can be computed as
\[
\frac{C_1 C_2}{C'_1 C'_2}
\frac{1}{\sqrt{2\pi}}
\int p(x\sqrt{C'_3}) y(x) e^{-x^2/2} dx,
\]
where
\[
y(x)
=
\sum_{0\le k}
\left(\frac{C_2\sqrt{C'_3}}{C'_2\sqrt{C_3}}\right)^k
\frac{h_k(x) h_k(0)}{k!}.
\]
Using the three-variable generating function, we find:
\[
e^{-x^2/2} y(x)
=
e^{-x^2/2} \sum_{k\ge 0} h_k(x) h_k(y) \frac{t^k}{k!}
=
C^{-1}
\exp(-(x/C)^2/2),
\]
where
\[
C := 
\frac{
\sqrt{(C'_2\sqrt{C_3})^2-(C_2\sqrt{C'_3})^2}
}{
C'_2 \sqrt{C_3}
},
\]
and we must satisfy the additional requirement
\[
|C_2 \sqrt{C'_3}|< |C'_2 \sqrt{C_3}|,
\]
precisely the additional hypothesis above.
But then
\[
q=
\frac{C_1 C_2}{C'_1 C'_2}
\frac{1}{\sqrt{2\pi}}
\int p(xC\sqrt{C'_3}) e^{-x^2/2} dx.
\]

Since
\[
\frac{C_1 C_2}{C'_1 C'_2} < 0,
\]
the problem reduces to the following.  For an integer $k>0$, how small can
we make the largest zero of a polynomial $p(x)$ of degree $2k+1$ subject to
the condition
\[
\int p(x) e^{-x^2/2} dx<0?
\]

Let $x^{(k)}_0$ be the smallest zero of $h_k(x)$, and consider the polynomials
\begin{align}
p_1(x) &= h_k(x)\\
p_2(x) &= \frac{h_k(x)}{x-x^{(k)}_0}.
\end{align}
By Christoffel-Darboux, we compute
\[
p_2(x) = h_{k-1}(x^{(k)}_0)^{-1} \sum_{0\le j<k} \frac{h_j(x) h_j(x^{(k)}_0)}{j!}
\]
We thus have the following integrals:
\begin{align}
\frac{1}{\sqrt{2\pi}} \int_{-\infty}^\infty p_1(x) p_2(x) e^{-x^2/2}dx
&=
0\\
\frac{1}{\sqrt{2\pi}} \int_{-\infty}^\infty p_2(x)^2 e^{-x^2/2} dx
&=
h_{k-1}(x^{(k)}_0)^{-2}
\sum_{0\le j<k} \frac{h_j(x^{(k)}_0)^2}{j!} > 0.
\end{align}
Now, consider the polynomial
\[
p(x) = p_1(x) p_2(x) - \epsilon p_2(x)^2 + \epsilon^2
     = (x-x^{(k)}_0-\epsilon) p_2(x)^2 + \epsilon^2
\]
for $\epsilon>0$ small.  This certainly satisfies the negative integral
condition; on the other hand, its only zero is that near $x^{(k)}_0$.

Thus when $l=2k$, we obtain a bound of the form:
\[
\limsup_{m\to\infty}
\frac{d_m-Lg(t'_0)m}{\sqrt{m}}
\le
C\sqrt{C'_3} x^{(k)}_0.
\]
Since
\[
x^{(k)}_0 = -2\sqrt{k} + O(k^{-1/6}),
\]
we are done.
\end{proof}

\begin{rems}
The additional assumption is satisfied in all of the applications above;
this is a trivial calculation for all but the lattice cases, in which it
follows from the transformation law.
\end{rems}

\begin{rems}
The fact that taking $2k$ coefficients gives an improvement proportional to
$2\sqrt{km}$ suggests that to obtain a first-order improvement, we will
need to let $k$ grow linearly with $m$.
\end{rems}

\begin{rems}
The involvement of the polynomial $h_k(x)^2/(x-x^{(k)}_0)$ is somewhat
reminiscent of the approach of Levenshtein \cite{LevenshteinVI:1998} to the
universa bounds, in which different orthogonal polynomials occur.  There
the polynomial giving the bound is
\[
(p_k(x)p_{k-1}(y)-p_k(y)p_{k-1}(x))^2/(x-y),
\]
which reduces to our polynomial when $p_k=h_k$ and
$y=x^{(k)}_0$.  While in those bounds, one lets $k$ grow linearly with $n$,
one can also consider finite $k$; in that case, the improvement is
$\Omega(\sqrt{n})$, as in our case.
\end{rems}

For the shadow cases, roughly the same argument applies; for instance, in
the self-dual binary code case, we take relations starting with
$m=n/12+(\log n)^2/6$.  Positivity near $j=m$ is immediate (since
multiplying by a positive power series leaves positive initial coefficients
positive); in the remaining region, we find that restricting the integral
to an interval $|x|\le n^{-\epsilon}$ and replacing $(1+t)^{(\log n)^2}$
by $(1+|t|)^{(\log n)^2}$ gives negligible relative error.  The argument
then proceeds as before.

For the quantum cases, the difficulty is in choosing the relation.
Basically, one defines $c^{(k)}_0$ and $d^{(k)}_0$ as above, and considers
a linear combination of
\[
m^{k/2} c^{(k)}_0\text{ and }
m^{(k+1)/2} ((q+1)c^{(k)}_0-(q-1)d^{(k)}_0).
\]
The first set of relations has no effect on the coefficients of $[t^j]
C(t)-D(t)$, so as above, we essentially obtain an arbitrary polynomial
here.  Similarly, they have a lower-order effect on the coefficient
of $[t^0] C(t)$; we thus end up with the same constraints on this polynomial
as above.  On the other hand, near $j/m=Lg(t'_0)$, the relations have
the same order behavior; in this neighborhood, we may thus choose an arbitrary
polynomial of degree $k+1$, so have no additional constraints.


\end{document}